\documentclass[graybox]{svmult}

\usepackage{mathptmx}       
\usepackage{helvet}         
\usepackage{courier}        
\usepackage{type1cm}        
\usepackage{makeidx}         
\usepackage{graphicx}        
\usepackage{multicol}        
\usepackage[bottom]{footmisc}

\usepackage{url}

\usepackage[english]{babel}     
\usepackage[utf8]{inputenc}  

\usepackage{geometry}
\usepackage{listings}
\usepackage{mathrsfs}
\usepackage{color}
\usepackage{amssymb}
\usepackage{afterpage}
\usepackage[T1]{fontenc}
\usepackage{amsfonts}
\usepackage{mathtools}
\usepackage{amsmath}
\usepackage{verbatim}
\usepackage{bm}
\usepackage{bbm}
\usepackage{enumerate}
\usepackage{stmaryrd}
\usepackage{tabu}
\usepackage{faktor}
\usepackage{algorithm,algcompatible}

\usepackage{subcaption}
\DeclareCaptionSubType*{algorithm}

\DeclareCaptionLabelFormat{alglabel}{Algorithm~#2}
\DeclareCaptionLabelFormat{algsublabel}{Algorithm~#2}

\graphicspath{{img/}}

\usepackage{moreverb} 

\newcommand{\rulesep}{\unskip\ \vrule\ }

\algnewcommand\algorithmicforeach{\textbf{for each}}
\algdef{S}[FOR]{FOREACH}[1]{\algorithmicforeach\ #1\ \algorithmicdo}

\newcommand*\mcupinn[2]{\vcenter{\hbox{$\mathsurround=0pt
  \ifx\displaystyle#1\textstyle\else#1\fi\bigcup$}}}

\newcommand*\mcapinn[2]{\vcenter{\hbox{$\mathsurround=0pt
  \ifx\displaystyle#1\textstyle\else#1\fi\bigcap$}}}
\DeclarePairedDelimiter{\abs}{\lvert}{\rvert}
\DeclarePairedDelimiter{\norm}{\lVert}{\rVert}

\DeclarePairedDelimiter{\bra}{\lbrace}{\rbrace}

\DeclarePairedDelimiter{\prodscal}{\langle}{\rangle}

\DeclareMathOperator*{\argmax}{arg\,max}
\DeclareMathOperator*{\expval}{\mathbb{E}}

\newcommand{\R}{{\mathbb R}}

\def \utruth {u_{\delta}}
\def \urb {u_N}
\def \rbspace {\mathbb{V}_N}
\def \trainset {\Xi_t}
\def \truthspace {\mathbb{V}_\delta}
\def \VV {\mathbb{V}}
\def \RR {\mathbb{R}}
\def \bmu {y}

\makeindex

\immediate\write18{texcount -inc -incbib -sum main.tex > .wordcount}

\begin{document}

\title*{Weighted reduced order methods for parametrized partial differential equations with random inputs}
\titlerunning{Weighted ROMs for parametrized PDEs with random inputs}
\author{Luca Venturi \and Davide Torlo \and Francesco Ballarin \and Gianluigi Rozza}
\authorrunning{L. Venturi \and D. Torlo \and F. Ballarin \and G. Rozza}
\institute{
Luca Venturi \at
mathLab, Mathematics Area, SISSA, Trieste, Italy \\
Current address: Courant Institute of Mathematical Sciences, New York University, U.S. \\ 
\email{venturi@cims.nyu.edu}
\and
Davide Torlo \at
mathLab, Mathematics Area, SISSA, Trieste, Italy \\
Current address: Institut f\"ur Mathematik. Universit\"at Z\"urich, Switzerland \\
\email{davide.torlo@math.uzh.ch}
\and
Francesco Ballarin \at
mathLab, Mathematics Area, SISSA, Trieste, Italy \\
\email{francesco.ballarin@sissa.it}
\and
Gianluigi Rozza \at
mathLab, Mathematics Area, SISSA, Trieste, Italy \\
\email{gianluigi.rozza@sissa.it}
}

\maketitle

\abstract{In this manuscript we discuss weighted reduced order methods for stochastic partial differential equations. Random inputs (such as forcing terms, equation coefficients, boundary conditions) are considered as parameters of the equations. We take advantage of the resulting parametrized formulation to propose an efficient reduced order model; we also profit by the underlying stochastic assumption in the definition of suitable weights to drive to reduction process. Two viable strategies are discussed, namely the weighted reduced basis method and the weighted proper orthogonal decomposition method. A numerical example on a parametrized elasticity problem is shown.}

\section{Introduction}
Several problems in applied sciences and engineering can be modeled by partial differential equations (PDEs). Complex models may be characterized by several coefficients, which could possibly be calibrated from experiments or measurements; as measurements of such quantities are often affected, in practice, by noise and uncertainty, one is usually interested in solving \emph{parametrized} PDEs, in which the coefficients are allowed to vary under a supposed probability distribution. As no analytical solution is usually available, numerical approximation by discretization methods (e.g. finite element) are customarily employed. However, it is often unaffordable (in terms of CPU times) to run a new simulation for each parametric scenario, especially when interested in the evaluation of some statistics of the solution. For this reason, several groups in the computational science and engineering community advocate the use of reduced order models.

Model reduction comes in different forms, often depending on the specific needs and applications. For instance, a dimensionality reduction of the physical space can be carried out to speed up computations; e.g. in fluid dynamics, engineers may occasionally be willing to replace a detailed model of a three-dimensional network of tubes by a simplified one-dimensional network, if local three-dimensional features are not relevant for the phenomenon at hand. The simplification can be carried out at the level of mathematical model as well, such as a potential flow description in fluid dynamics. This work will focus, instead, on \emph{projection-based} reduction techniques, which are based on the idea that parametrized solutions often belong to a lower dimensional manifold of the solution space. In such a case, few basis functions (obtained from few solutions of the detailed physical problem) should be enough to span the lower dimensional manifold.

Two widely used techniques for basis functions generation are going to be discussed in this work, namely the reduced basis method \cite{rozza:book,reduced_50,rozza_nguyen,heat_tranfer_rozza_huynh} and the proper orthogonal decomposition approach \cite{rozza:book,ito}. While originally proposed for deterministic problems, such methods have been extended as well to a stochastic framework by properly weighting the reduced basis construction phase according to the underlying probability distribution \cite{peng,peng_art,peng_A}. The resulting reduced order methods (ROMs) have been named \emph{weighted reduced basis} method \cite{peng, peng_art,peng_A,torlo2017stabilized,spannring2017weighted,spannring2018phd} and \emph{weighted proper orthogonal decomposition} \cite{VenturiBallarinRozza2018}. The main goal of this chapter is to summarize and compare these two weighted approaches by means of a numerical test case. We do mention that further ROM techniques to tackle uncertainty are available, e.g. \cite{hestaven:collocation,gunzburger,Traian3,nobile:collocation}, but a more extensive comparison is beyond the scope of this work.

The outline of the chapter is as follows: the mathematical model, governed by a parametrized PDE, is outlined in Section \ref{sec:fe}, as well as its finite element discretization. Two weighted ROMs are proposed in Section \ref{sec:wrom}, which extend deterministic reduced basis and proper orthogonal decomposition methods. Numerical results on a elasticity problem test case are provided in Section \ref{sec:results}. Conclusions follow in Section \ref{sec:conclusion}.

\section{Parametrized formulation and high fidelity approximation}
\label{sec:fe}

Let $(\Omega, \mathcal{F}, P)$ be a complete probability space, where $\Omega$ denotes the set of possible outcomes, $\mathcal{F}$ is the $\sigma$-algebra of events and $P$ is the probability measure. Let $y = y(\omega)$, $y: (\Omega, \mathcal{F}) \to (\Gamma, \mathcal{B})$, where $\Gamma \subset \mathbb{R}^K$ is a compact set, $\mathcal{B}$ is the Borel measure, and $y(\omega) = (y_1(\omega), \hdots, y_K(\omega))$ is a random vector which components are independent absolutely continuous random variables. Denote by $\rho: \mathbb{R}^K \to \mathbb{R}$ the probability density function of $y$. In the following, we will identify the realization of $y$ as the parameter.

Furthermore, let $D \subset \mathbb{R}^d$, $d = 1, 2, 3$, be a bounded domain. We assume that the engineering problem at hand is modeled by a parametrized elliptic partial differential over $D$, as follows:
\begin{align*}
&\text{find }u: \Gamma \to \mathbb{V}\text{ such that}\\
&a(u(y(\omega)), v; y(\omega)) = F(v; y(\omega)) \quad \forall v \in \mathbb{V},
\end{align*}
for a.e. $\omega\in \Omega$.

In the following we will assume $\mathbb{V}$ to be a subset of $H^1(D)$; we will also assume that the bilinear form $a(\cdot, \cdot; y): \mathbb{V} \times \mathbb{V} \to \mathbb{R}$ is coercive and continuous, as well as the linear functional $F(\cdot; y) \in \mathbb{V}'$, in order to guarantee the well-posedness of the problem \cite{quarteroni_valli}. Specific expressions for $a(\cdot, \cdot; y)$ and $F(\cdot; y)$ will be provided in the numerical experiments in Section \ref{sec:results}.

In general, due to the lack of analytical solutions, we will resort to a numerical approximation based on a Galerkin approach.
Given a finite dimensional approximation space $\mathbb{V}_\delta \subseteq \mathbb{V}$ (in the following, a finite element space \cite{quarteroni_valli}), characterized by $\mathrm{dim}(\mathbb{V}_\delta) = N_\delta < \infty$, we consider the approximate problem:
\begin{align}
&\text{find }u_{\delta}: \Gamma \to \mathbb{V}_\delta\text{ such that}\notag\\
&a(u_{\delta}(y(\omega)),v;y(\omega)) = F(v;y(\omega)) \qquad \text{for all } v\in\mathbb{V}_\delta,
\label{stochastic_elliptic_pde_FE}
\end{align}
for a.e. $\omega \in \Omega$. We refer to problem \eqref{stochastic_elliptic_pde_FE} as the \emph{high fidelity} (or \emph{truth}) \emph{problem} and to $u_{\delta}$ as the \emph{high fidelity solution}.

We are interested in the computation of statistics of the solution, e.g. its expectation $\expval[u]$, or statistics of some output of interest $s: \mathbb{V} \to \mathbb{R}$, e.g. $\expval[s(u)]$, by means of Monte-Carlo method. Since the repeated evaluation of \eqref{stochastic_elliptic_pde_FE} for different realizations $y \in \Gamma$ is usually computationally expensive for moderately large $N_\delta$ (i.e., high accuracy of $u_{\delta}(y) \simeq u(y)$), we aim at replacing the \emph{high fidelity} solution $u_{\delta}(y)$ with a cheaper, yet accurate, approximation. This approximation will be the result of a \emph{reduced order method} in which a Galerkin project is carried out over a further low dimensional space $\VV_N$ (rather than $\VV_\delta$).

\section{Weighted reduced order methods}
\label{sec:wrom}

In this Section we introduce the \emph{weighted} reduced order methods (ROMs), which extends deterministic ROMs by means of proper weighting associated to the underlying probability space. As usual in the ROM literature, we ensure an offline-online computational strategy; during the offline stage the reduced basis space $\VV_N$ is constructed once and for all, while during the online phase (to be repeated for any possible realization of $y$) a Galerkin projection over $\VV_N$ is carried out to find the \emph{reduced order solution}. For the sake of brevity, in this Section we will only discuss in detail the offline phase, as it is the one in which a proper weighting is needed in order to account for stochastic formulation of the problem. The online phase, instead, can be carried out as in the deterministic setting; we refer to \cite{rozza:book,rozza_nguyen,heat_tranfer_rozza_huynh} for further details.

\subsection{Weighted Reduced Basis method}
Let us introduce first the weighted Reduced Basis method, originally introduced in \cite{peng,peng_art,peng_A} as an extension of the (deterministic) Reduced Basis (RB) method \cite{rozza:book}.

\begin{algorithm}%
\captionsetup{labelformat=alglabel}
\caption{Greedy algorithms}%
\label{algo:greedy}%
\begin{subalgorithm}{.49\textwidth}%
\captionsetup{labelformat=algsublabel}%
\caption{(Deterministic) Greedy Algorithm}%
\label{algo:dgreedy}%
 \hspace*{\algorithmicindent} \textbf{Input:} parameter domain $\Gamma$, tolerance $\varepsilon^{tol}$ and $N_{max}$.\\
 \hspace*{\algorithmicindent} \textbf{Output:} reduced space $\rbspace$.
\begin{algorithmic}[1]%
{\STATE Sample $\trainset \subset \Gamma$;
 \STATE Define $\mathbb{V}_0 = \emptyset$;
 \STATE Pick arbitrary $y^1 \in \trainset$;
 \FOR{$N = 1, \hdots, N_{max}$}
 \STATE Solve \eqref{stochastic_elliptic_pde_FE} for $y = y^N$ to compute $u(y^N)$;
 \STATE Update $\VV_{N} = \VV_{N-1} \bigoplus \text{span}\bra{ \utruth (y^N)}$;
 \STATE Compute $y^{N+1}=\argmax_{y\in \trainset} \eta_N(y)$;
 \IF{$\eta_N(y^{N+1})\leq\varepsilon^{tol}$}
 \STATE \textbf{break}
 \ENDIF
\ENDFOR
}
 \end{algorithmic}
\end{subalgorithm}%
\rulesep
\begin{subalgorithm}{.49\textwidth}%
\captionsetup{labelformat=algsublabel}%
\caption{Weighted Greedy Algorithm}%
\label{algo:wgreedy}%
 \hspace*{\algorithmicindent} \textbf{Input:} parameter domain $\Gamma$, tolerance $\varepsilon^{tol}$ and $N_{max}$.\\
 \hspace*{\algorithmicindent} \textbf{Output:} reduced space $\rbspace$.
\begin{algorithmic}[1]%
{\STATE \emph{Properly} sample $\trainset \subset \Gamma$;
 \STATE Define $\mathbb{V}_0 = \emptyset$;
 \STATE Pick arbitrary $y^1 \in \trainset$;
 \FOR{$N = 1, \hdots, N_{max}$}
 \STATE Solve \eqref{stochastic_elliptic_pde_FE} for $y = y^N$ to compute $u(y^N)$;
 \STATE Update $\VV_{N} = \VV_{N-1} \bigoplus \text{span}\bra{ \utruth (y^N)}$;
 \STATE Compute $y^{N+1}=\argmax_{y\in \trainset} \eta_N^{w}(y)$;
 \IF{$\eta_N^{w}(y^{N+1})\leq\varepsilon^{tol}$}
 \STATE \textbf{break}
 \ENDIF
\ENDFOR
}
 \end{algorithmic}
\end{subalgorithm}
\end{algorithm}

In the deterministic setting, the construction of the reduced space $\VV_N$ is based on the so-called Greedy algorithm, summarized in Algorithm \ref{algo:dgreedy}. The aim of the Greedy algorithm is to build a reduced basis space $\rbspace$ spanned by some of the \textit{truth} snapshots $\bra{ \utruth (y^i) }_{i=1}^{n_t} $, where $\trainset = \bra{y^i}_{i=1}^{n_t} \subset \Gamma$ is a training set of parameters (line 1). After an initialization to setup the empty reduced space (line 2) and the first value $y^1$ (line 3), the algorithm proceeds iteratively (lines 4-11). At each iteration $N$, the solution corresponding to the current parameter $y^N$ is computed (line 5) and used to enrich the reduced space (line 6). The parameter $y^{N+1}$ to be used at the next iteration is then automatically chosen (in a greedy fashion) as the worst approximated solution on the whole set of parameters $\trainset$ (line 7). The simplest implementation of such selection criterion would require the computation of the error $e_N(y)= || \utruth(y) - \urb (y)||_{\mathbb{V}}$ for each $y \in \trainset$, and the corresponding computation of
\begin{equation}
\arg\max_{\bmu\in \trainset} ||\utruth(\bmu)-\urb (\bmu)||_{\VV}.
\label{eq:dargmax}
\end{equation}
A more efficient alternative is to use an error estimator $\eta_N(y)$ that bounds from above the error, i.e. such that $$e_N(y) \leq \eta_N(y),$$ and that does not require the evaluation of the \textit{truth} solution $\utruth (y)$. This is readily possible for parametrized elliptic problems as in Section \ref{sec:fe}, under mild assumptions \cite{rozza:book}. Thus, the computation of
\begin{equation*}
\arg\max_{y \in \trainset} \eta_N(y)
\end{equation*}
is used in place of \eqref{eq:dargmax}.

\par In a stochastic context, we want to give different importance to different parameters (i.e., different realizations), reflecting the probability distribution. This can be easily pursued by using a different norm $||\cdot||_w$ in \eqref{eq:dargmax}, i.e. by optimizing
\begin{equation}
\arg\max_{\bmu\in \trainset} ||\utruth(\bmu)-\urb (\bmu)||_{w},
\label{eq:wargmax}
\end{equation}
where
\begin{equation}
\label{stochastic_norm_weighted_space}
||u(\bmu)||_w=w(\bmu)||u(\bmu)||_{\VV} \quad \forall u \in \truthspace,\forall \bmu \in \Gamma,
\end{equation} 
and $w : \Gamma \to \mathbb{R}^+$ is a suitable weight function. The choice of the weight function $w(\bmu)$ may be driven by the desire of minimizing the expectation of the square norm error of the RB approximation, i.e.
\begin{equation}
\label{square_RB_w_error}
\begin{split}
\mathbb{E}[||\utruth -\urb||_{\VV}^2]=&\int_\Omega ||\utruth(\bmu(\omega))-\urb(\bmu(\omega)))||_{\VV}^2 dP(\omega)=\\
=&\int_{\Gamma} ||\utruth(\bmu)-\urb(\bmu)||_{\VV}^2  \rho(\bmu) d\bmu,
\end{split}
\end{equation}
that we can bound with
\begin{equation}
\label{RB_w_error_bounding}
\mathbb{E}\left[||\utruth-\urb||_{\VV}^2\right]\leq \int_{\Gamma} \eta_N(\bmu)^2 \rho(\bmu) d\bmu.
\end{equation}
This suggests to define the weight function as $w(\bmu)=\sqrt{\rho(\bmu)}$ as proposed in \cite{peng_art}.

Other choices of the weight function $w$ are possible. For instance, let us suppose that we are interested in the accurate approximation (on average) of a quantity $s :\VV \to \R$, which depends on the solution of the parametric PDE $u$, namely in minimizing
\begin{equation}
    \expval[|s(\utruth)-s(\urb)|]=\int_{\Gamma} |s(\utruth(y))-s(\urb(y)| \rho(y) dy.
\end{equation}
In this case, it is natural to choose as weight function $w(y)=\rho(y)$. 

Regardless of the choice of $w$, the first modification of the weighted approach lies in using
the weighted error estimator
\begin{equation}
\eta_N^w(\bmu) = w(y) \eta_N(y)
\end{equation}
on line 7 of Algorithm \ref{algo:wgreedy}.

\par
Furthermore, the choice of training set $\Xi_{t}$ is usually crucial in the stochastic context (line 1 of Algorithm \ref{algo:wgreedy}). Two choices that we will compare in the numerical examples are (1) sampling from the uniform distribution $\mathcal{U}(\Gamma)$, as it is usually done in the deterministic case, and (2) sampling from the distribution $\rho$.

\subsection{Weighted Proper Orthogonal Decomposition}
\label{weighted_POD}

Let us introduce next the weighted POD method, originally introduced in \cite{VenturiBallarinRozza2018} as an extension of the deterministic POD--Galerkin approach \cite{rozza:book}. For instance, the weighted POD may be used when an error estimator $\eta_N(y)$ is not readily available for the weighted RB approach.

\begin{algorithm}%
\captionsetup{labelformat=alglabel}
\caption{POD algorithms}%
\label{algo:pod}%
\begin{subalgorithm}{.49\textwidth}%
\captionsetup{labelformat=algsublabel}%
\caption{(Deterministic) POD algorithm}
\label{algo:dpod}
 \hspace*{\algorithmicindent} \textbf{Input:} parameter domain $\Gamma$, tolerance $\varepsilon^{tol}$ and $N_{max}$.\\
 \hspace*{\algorithmicindent} \textbf{Output:} reduced space $\rbspace$.
\begin{algorithmic}[1]%
\STATE Sample $\trainset \subset \Gamma$;
\FOREACH{$y \in \trainset$}
 \STATE Solve \eqref{stochastic_elliptic_pde_FE} to compute $u(y)$;
\ENDFOR
\STATE Assemble the matrix $\hat{C}_{ij} = \prodscal{\varphi_i,\varphi_j}_\mathbb{V}$;
\STATE Compute its eigenvalues $\lambda_k$, and corresponding eigenvectors $\psi_k$, $k = 1, \hdots, n_t$;
\STATE Find the minimum $N \in \lbrace  1, \hdots, N_{max} \rbrace$ such that $E_{N} > 1 - \varepsilon^{tol}$;
\STATE Define $\mathbb{V}_N = \mathrm{span}\bra{\xi^1,\dots,\xi^N}$, being $\xi^i = \sum_{j=1}^{n_t} \psi^i_j \varphi_j$.
\end{algorithmic}
\end{subalgorithm}%
\rulesep
\begin{subalgorithm}{.49\textwidth}%
\captionsetup{labelformat=algsublabel}%
\caption{Weighted POD algorithm}
\label{algo:wpod}%
 \hspace*{\algorithmicindent} \textbf{Input:} parameter domain $\Gamma$, tolerance $\varepsilon^{tol}$ and $N_{max}$.\\
 \hspace*{\algorithmicindent} \textbf{Output:} reduced space $\rbspace$.
\begin{algorithmic}[1]%
\STATE \emph{Properly} sample $\trainset \subset \Gamma$;
\FOREACH{$y \in \trainset$}
 \STATE Solve \eqref{stochastic_elliptic_pde_FE} to compute $u(y)$;
\ENDFOR
\STATE Assemble the \emph{weighted} matrix $\hat{C}^w_{ij} = w_i\prodscal{\varphi_i,\varphi_j}_\mathbb{V}$;
\STATE Compute its eigenvalues $\lambda_k$, and corresponding eigenvectors $\psi_k$, $k = 1, \hdots, n_t$;
\STATE Find the minimum $N \in \lbrace 1, \hdots, N_{max} \rbrace$ such that $E_{N} > 1 - \varepsilon^{tol}$;
\STATE Define $\mathbb{V}_N = \mathrm{span}\bra{\xi^1,\dots,\xi^N}$, being $\xi^i = \sum_{j=1}^{n_t} \psi^i_j \varphi_j$.
\end{algorithmic}
\end{subalgorithm}%
\end{algorithm}

In its deterministic version, the POD method aims to minimize the mean square error:
\begin{equation}
\label{pod_deterministic_l2_error}
    \int_\Gamma \norm{u(y) - u_N(y)}_\mathbb{V}^2\,dy.
\end{equation}
over all possible reduced spaces of dimension $N$. From a practical point of view, a training set $\trainset$ of cardinality $n_t$ is introduced (line 1 of Algorithm \ref{algo:dpod}), and high fidelity solutions are computed for each $y \in \trainset$ (lines 2-4).
This enables to consider a discretized version of \eqref{pod_deterministic_l2_error}:
\begin{equation}
\frac{1}{n_t}\sum_{y\in \Xi_t} \norm{u(y) - u_N(y)}_\mathbb{V}^2
\label{pod_deterministic_l2_error_discretized}
\end{equation}
where $\varphi_i = u(y^i)$, $i=1,\dots,n_t = \abs{ \Xi_t}$. The reduced order subspace $\mathbb{V}_N$ is therefore defined as the space spanned by the $N$ leading eigenvectors of the operator \cite{rozza:book}
$$
C : v\in\mathbb{V} \to \sum_{i=1}^{n_t} \prodscal{v,\varphi_i}_\mathbb{V}\cdot \varphi_i.
$$
In practice, a basis of $\mathbb{V}_N$ can be found by computing the leading eigenvectors of the matrix $\hat{C} = (\prodscal{\varphi_i,\varphi_j}_\mathbb{V})_{1\leq i,j\leq n_t}$ (lines 5-6). Denoting by ($\lambda_k$, $\psi_k$) the eigenvalue-eigenvector pairs, the optimal subspace $\VV_N$ is obtained as
$$\mathbb{V}_N = \mathrm{span}\bra{\xi^1,\dots,\xi^N},$$
being $\xi^i = \sum_{j=1}^{n_t} \psi^i_j \varphi_j$ (line 8). Practical stopping criterion involves either a maximum reduced space dimension $N_{max}$, or the selection of the minimum dimension $N$ such that the so called retained energy
$$E_N = \frac{\sum_{k=1}^{N} \lambda_k}{\sum_{k=1}^{n_t} \lambda_k}$$
is greater than a prescribed tolerance (line 7).

To accommodate a stochastic setting, we observe that \eqref{pod_deterministic_l2_error_discretized} involves a uniform Monte-Carlo approximation of \eqref{pod_deterministic_l2_error}. 
In the stochastic case, one would like to find the reduced order sub-space $\mathbb{V}_N$ which minimizes the mean square error
\begin{equation}
\mathbb{E}[\norm{u -u_N }_\mathbb{V}^2] = \int_\Gamma \norm{u(y(\omega))-u_N(y(\omega))}_\mathbb{V}^2\,dP(\omega).
\label{pod_stochastic_l2_error}
\end{equation}
Nonetheless, any quadrature rule can be used for that purpose. In this framework, let $\Xi_t \subset \Gamma$ and $w_i$, $i = 1, \hdots, n_t$ be, respectively, the set of nodes and weights of a chosen quadrature rule:
\begin{equation}
\label{pod_deterministic_l2_error_quad_rule}
    \int_\Gamma \norm{u(y(\omega)) - u_N(y(\omega))}_\mathbb{V}^2\,dP(\omega) \simeq \sum_{i=1}^{n_t} w_i \norm{u(y_i) - u_N(y_i)}_\mathbb{V}^2.
\end{equation}
Computationally this turns out to require the evaluation of the eigenvalues of a weighted (preconditioned) matrix $$C^w = (\omega_i\prodscal{\varphi_i,\varphi_j}_\mathbb{V})_{1\leq i,j\leq n_t} = W C \quad \text{for} \quad W = \mathrm{diag}(w_1,\dots,w_{n_t})$$
We observe that the matrix $C^w$ is not diagonalizable in the usual sense, but with respect to the scalar product induced by $C$; therefore this allows to recover $N$ orthogonal leading eigenvectors.

\begin{figure}
\centering
\begin{minipage}[c]{.47\textwidth}
\includegraphics[width=\textwidth,
keepaspectratio]{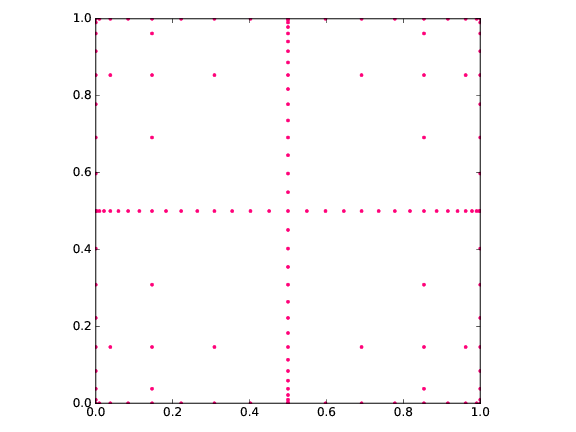}
\end{minipage}
\hspace{4mm}
\begin{minipage}[c]{.47\textwidth}
\includegraphics[width=\textwidth,
keepaspectratio]{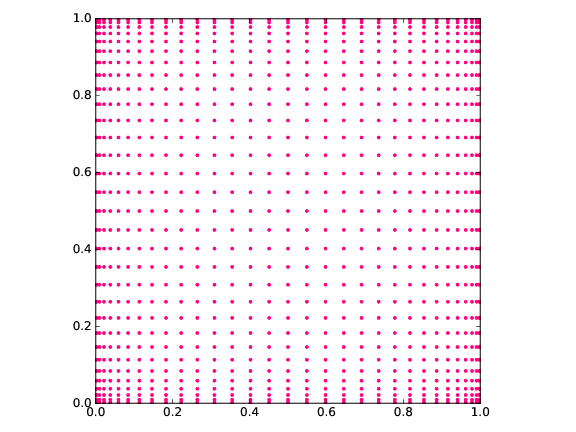}
\end{minipage}
\caption{Two dimensional grids based on nested Clenshaw-Curtis nodes of order $q=6$. The left one is based on a Smolyak rule ($145$ nodes), while the right one on a full tensor product rule ($1089$ nodes).}
\label{full_and_sparse_quadrature_rules_figure}
\end{figure}

We note that, once a quadrature rule has been selected, the weighted POD method prescribes both the weights (employed in line 5 of Algorithm \ref{algo:wpod}) and the training set (line 1). This is a difference to the weighted RB method, which only prescribes weights without an associated training set. Concerning the choice of the quadrature rule, the most straightforward way to do this is to use a Monte-Carlo approximation of \eqref{pod_stochastic_l2_error}, i.e. to choose $\Xi_t$ sampled from the distribution $\rho$ and $w_i \equiv 1$. For higher dimensional parameter spaces, a simple option is to perform a tensor product of a one-dimensional quadrature rules; however, in order to not be affected by the curse of dimensionality (which would entail very large training sets), we can instead adopt a Smolyak-type sparse quadrature rule \cite{novak:sparse_grids,holtz:sparse_grids,wasil:sparse_grids,gerstner:sparse_grids,novak:sparse_grids_quadrature} to significantly reduce the computational load (see e.g. \figurename~\ref{full_and_sparse_quadrature_rules_figure})

\section{Numerical results}
\label{sec:results}
In this Section we apply weighted reduced order methods to a linear elasticity problem.

\begin{figure}
\centering
\includegraphics[width=0.5\textwidth]{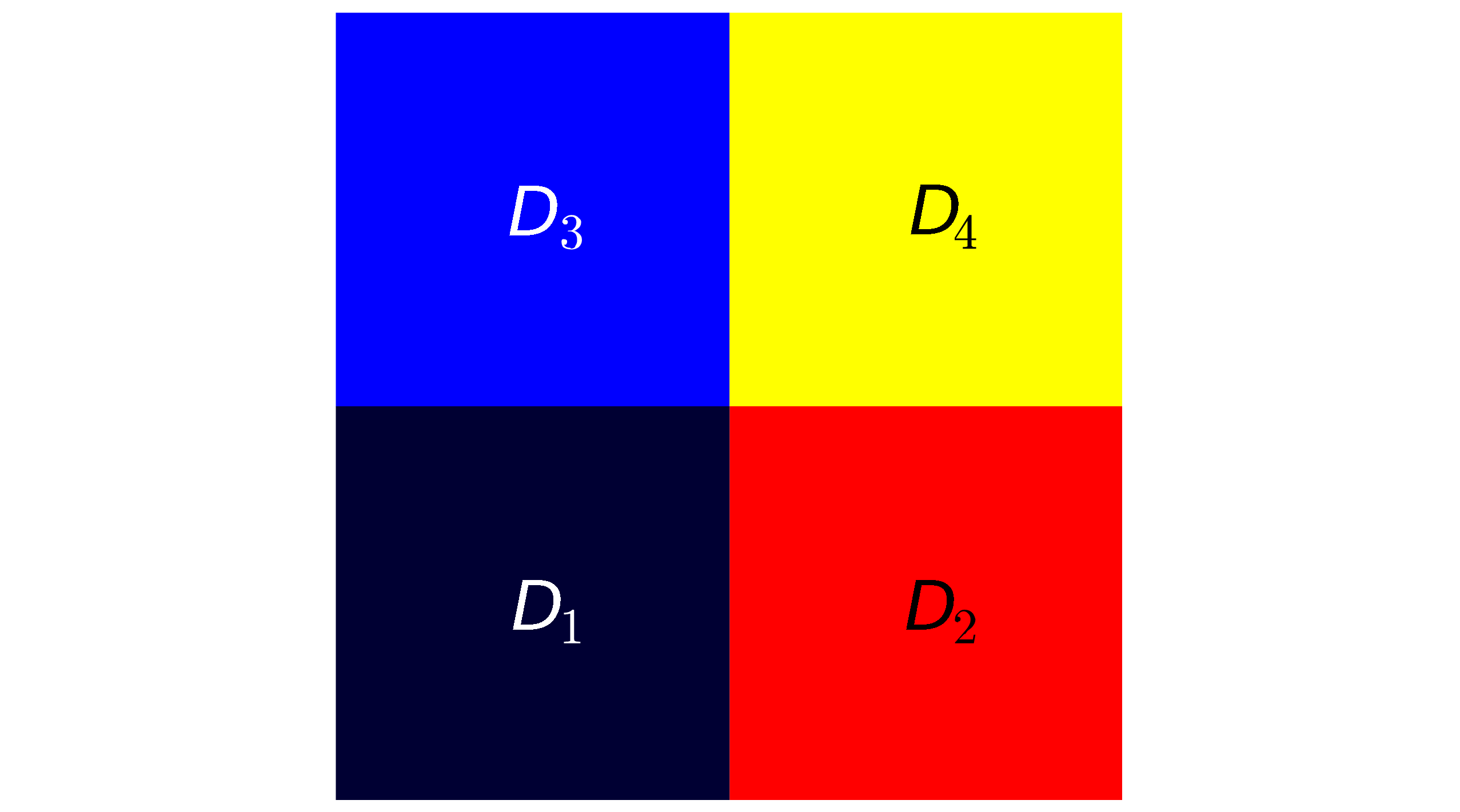}
\caption{Geometrical set-up of linear elasticity test case.}
\label{fig:domain}
\end{figure}

Let $D_1, \hdots, D_4$ be a partition of the domain $D = [0, 1]^2$, as in \figurename~\ref{fig:domain}. We consider a linear elasticity problems characterized by six parameters, as follows:
\begin{align*}
a(\mathbf{u},\mathbf{v}; y^1,\dots,y^4) & = \sum_{i=1}^4 y^i \int_{D_i} \{ \lambda \, (\nabla\cdot\mathbf{u})(\nabla\cdot\mathbf{u}) + 2\mu \, \mathbf{e}(\mathbf{u})\,\mathbf{:}\,\mathbf{e}(\mathbf{v}) \}\, dx, \\
f(\mathbf{v}; y^5, y^6) & = \sum_{i=1}^2 y^{i+4} \int_{(i-1)/2}^{i/2} v_2(1,x_2) \, dx_2,
\end{align*}
for every $\mathbf{u},\mathbf{v}\in\mathbb{V}$, where $\mathbb{V}=\{ \mathbf{v} \in H^1(D; \RR^2) \; : \; \mathbf{v}|_{[0,1] \times \{0,1\}} = 0 \}$, $\nabla \cdot \mathbf{u}$ denotes the divergence of the displacement $\mathbf{u}$, and $\mathbf{e}(\mathbf{u})$ is the symmetric part of the gradient of $\mathbf{u}$.
This corresponds to consider the equation of a linear elastic block, split in four parts, with Lam\'e constants rescaled by a parameter $y^i$, $i=1,\dots,4$. The unscaled Lam\'e constants $\lambda$ and $\mu$ are obtained for Young modulus equal to 1 and Poisson ratio set to 0.3. Boundary conditions correspond to later tractions equal to $y^5$ ($y^6$) on the bottom (top) half of the right side, and clamped on bottom and top sides.

The random variables $y^i$, $i=1,\dots,6$, are such that  
\begin{align*}
\frac{y^i-1}{2} & \sim \mathrm{Beta}(\alpha, \beta), \qquad \textrm{for $i=1,\dots,4$,} \\
\frac{y^i-2}{4} & \sim \mathrm{Beta}(\alpha, \beta), \qquad \textrm{for $i=5,6$.}
\end{align*}
We will consider two cases, namely $(\alpha, \beta) = (10, 10)$ and $(\alpha, \beta) = (75, 75)$, the latter resulting in a more peaked probability distribution.

\begin{figure}
\centering
\begin{minipage}[c]{.47\textwidth}
\includegraphics[width=\textwidth,
keepaspectratio]{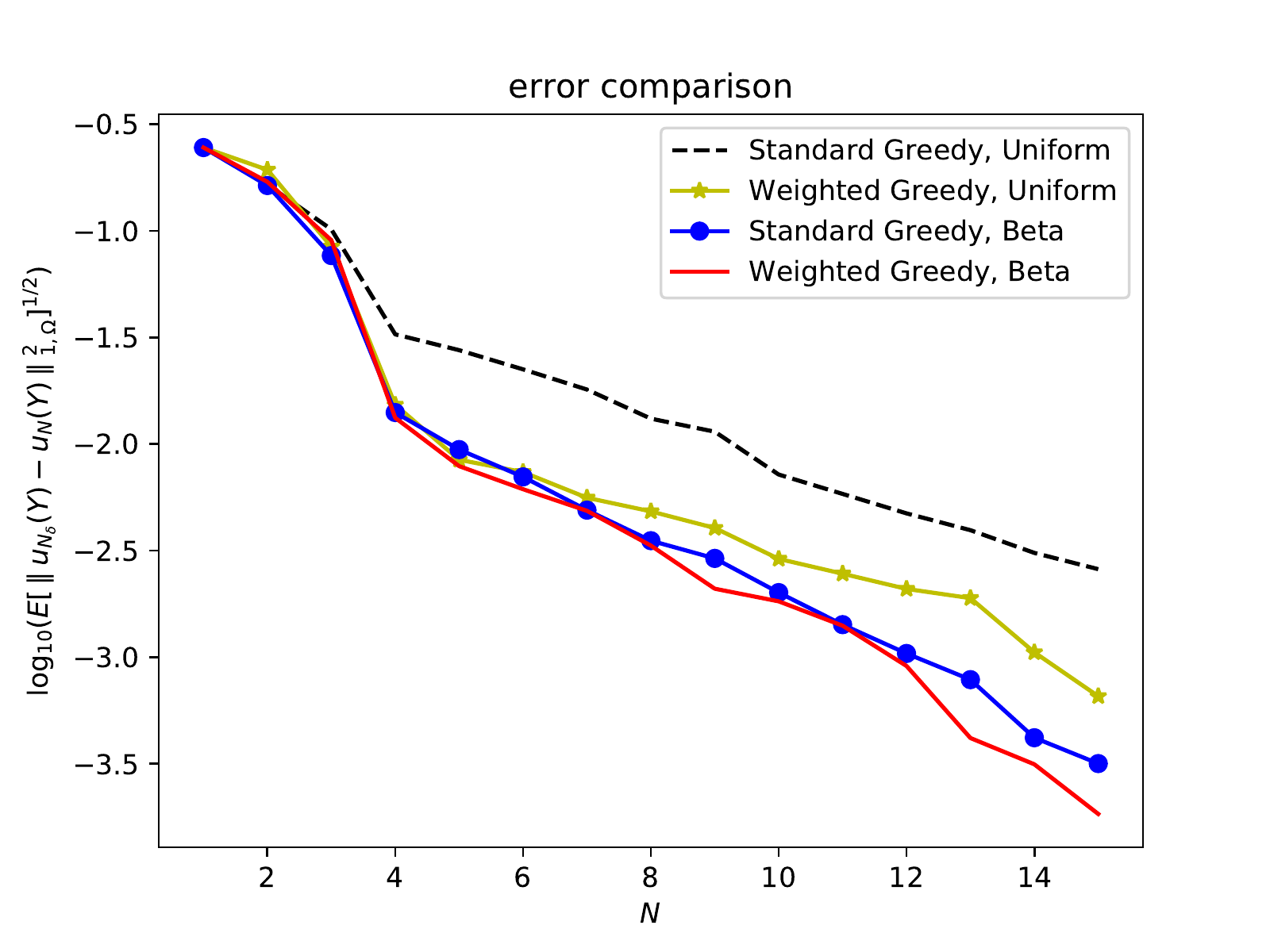}
\end{minipage}
\hspace{4mm}
\begin{minipage}[c]{.47\textwidth}
\includegraphics[width=\textwidth,
keepaspectratio]{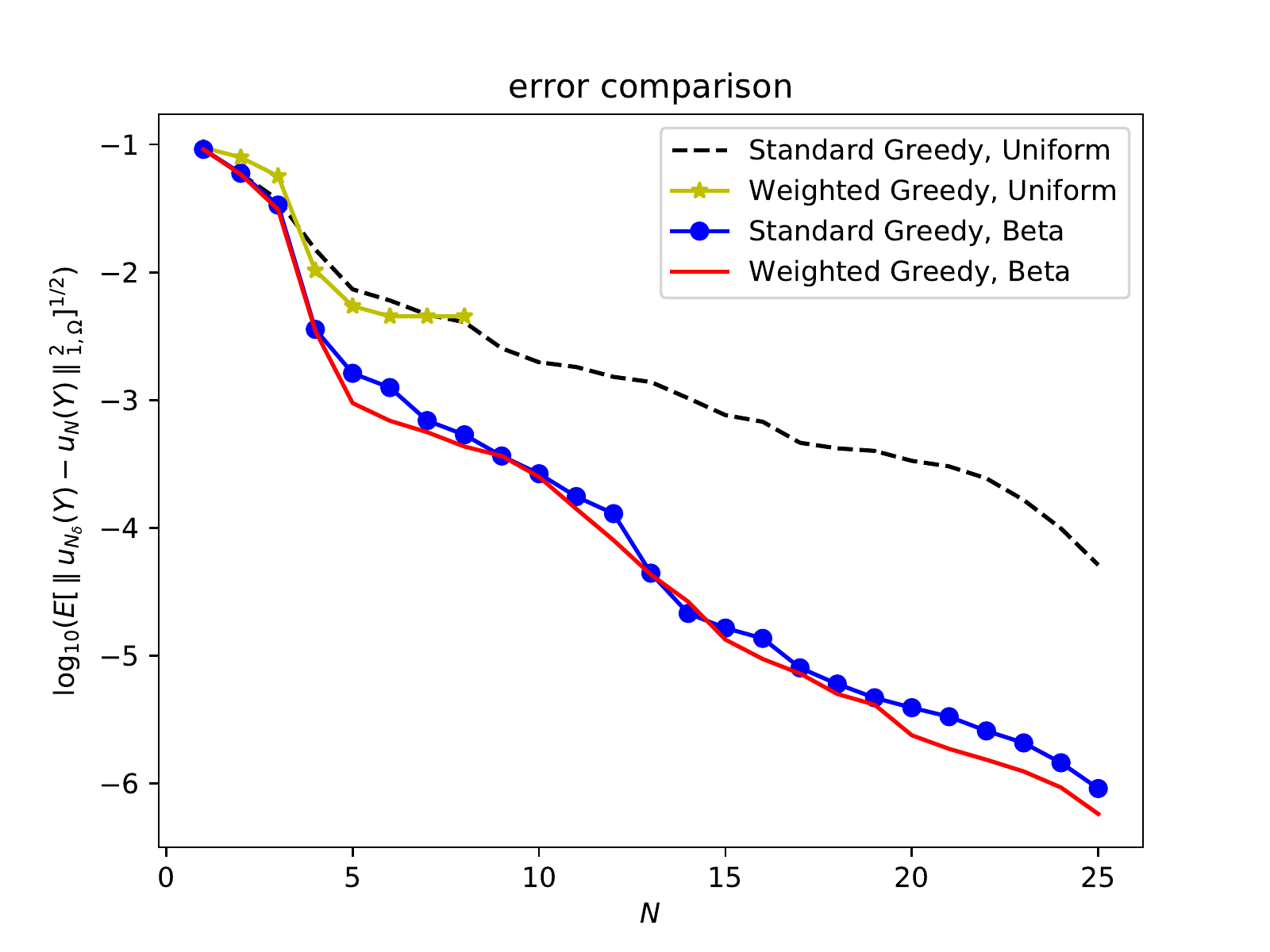}
\end{minipage}
\caption{Error analysis for standard and weighted RB algorithms, employing either a Uniform or Beta sampling for $\trainset$, for the cases $(\alpha, \beta)=(10, 10)$ (left) and $(\alpha, \beta)=(75, 75)$ (right). \label{figure:6E_rb}}
\end{figure}

Let us first comment the results of the weighted RB approach, summarized in \figurename~\ref{figure:6E_rb}, in which we plot the error $\mathbb{E}[\norm{u -u_N }_\mathbb{V}^2]$ for increasing reduced basis dimension $N$. We consider four different cases, which concern both the weighting (standard Greedy vs weighted Greedy) and the sampling (Uniform sampling vs sampling from the underlying Beta distribution) used during the Offline stage. \figurename~\ref{figure:6E_rb}(left), obtained for $(\alpha, \beta) = (10, 10)$, shows that, with an Uniform sampling, a proper weighting allows to improve the accuracy of the reduced order model (compare the dashed black line to the solid gold one). Further improvement can be obtained by sampling from the underlying Beta distribution (solid blue and solid red lines). The best accuracy is obtained if both the correct weighting and the appropriate sampling are employed (solid red line), resulting in an improvement of more than one order of magnitude compared to the standard Greedy algorithm employed in the deterministic setting. Such an improvement is even more prominent in \figurename~\ref{figure:6E_rb}(right), which corresponds to $(\alpha, \beta) = (75, 75)$; in this case, the best weighted approach outperforms the deterministic algorithm by two orders of magnitude. The importance of a proper sampling is even more remarkable in this case, as the weighted Greedy with Uniform sampling suffered a breakdown after very few iterations due to the improper exploration of the parameter space, which led to a singular reduced order system.

\begin{table}
\centering
{\renewcommand\arraystretch{1.2}
\begin{tabular}{|l|c|c|c|c|c|}
\hline
 & Standard & Uniform Monte-Carlo & Monte-Carlo & Gauss-Legendre & Sparse Gauss-Jacobi \\ \hline
$\abs{\Xi_t}$ & $500$ & $500$ & $500$ & $729$ & $389$
\\ \hline
\end{tabular}} 
\caption{Cardinality of the training set $\Xi_t$ employed in the various POD algorithms.\label{table:4EB_POD}}
\end{table}

\begin{figure}
\centering
\begin{minipage}[c]{.47\textwidth}
\includegraphics[width=\textwidth,
keepaspectratio]{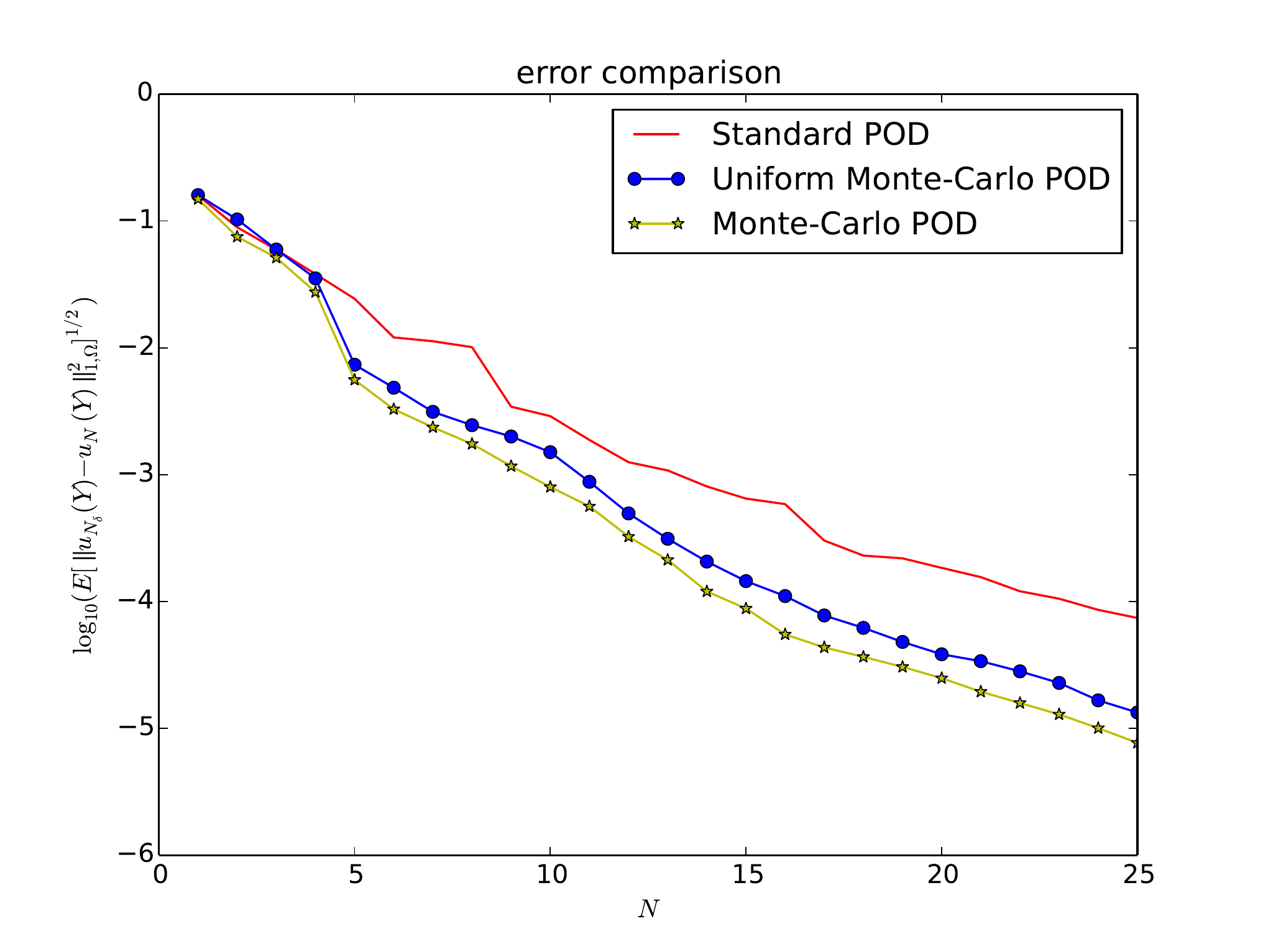}
\end{minipage}
\hspace{4mm}
\begin{minipage}[c]{.47\textwidth}
\includegraphics[width=\textwidth,
keepaspectratio]{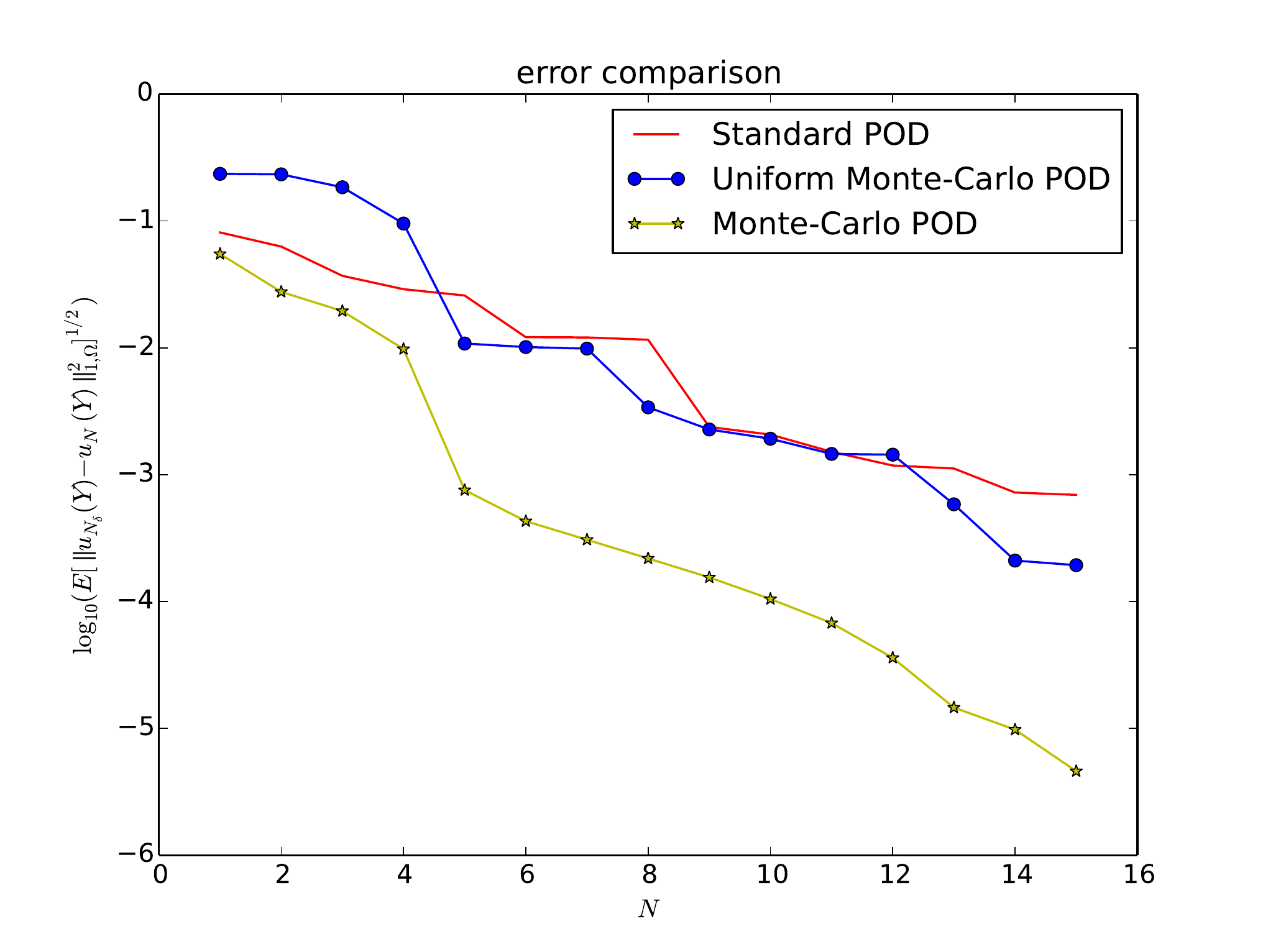}
\end{minipage}
\caption{Error analysis for standard, Uniform Monte-Carlo and Monte-Carlo POD algorithms, for the cases $(\alpha, \beta)=(10, 10)$ (left) and $(\alpha, \beta)=(75, 75)$ (right). \label{figure:6E_mcPOD}}
\end{figure}

We next carry out a similar task for different POD algorithms. \figurename~\ref{figure:6E_mcPOD} reports a comparison between standard POD (i.e. sampling from Uniform distribution, no weighting), Uniform Monte-Carlo POD (i.e. sampling from a Uniform distribution, and weighting according to the probability density function $\rho$) and Monte-Carlo POD (i.e. sampling from the underlying Beta distribution, without any weighting). In contrast to the previous comparison for RB methods, we note that the fourth case (i.e. sampling from Beta and weighting according to $\rho$) is omitted, as it does not admits an interpretation in terms of the quadrature formula \eqref{pod_deterministic_l2_error_quad_rule}. The cardinality of the training set (i.e., the number of quadrature points) is set to 500 for all three cases (see \tablename~\ref{table:4EB_POD}). The results in \figurename~\ref{figure:6E_mcPOD} show that Monte-Carlo POD guarantees the best results, with improvements of one (left) and two (right) orders of magnitude, respectively, with respect to the standard (deterministic) POD in terms of accuracy. Results for $(\alpha, \beta) = (75, 75)$ emphasize the importance of a correct sampling (and thus, correct weighting as induced by \eqref{pod_deterministic_l2_error_quad_rule}) as well: indeed, the Uniform Monte-Carlo POD provides marginal improvements over the standard POD in the case of highly concentrated distribution.

\begin{figure}
\centering
\begin{minipage}[c]{.47\textwidth}
\includegraphics[width=\textwidth,
keepaspectratio]{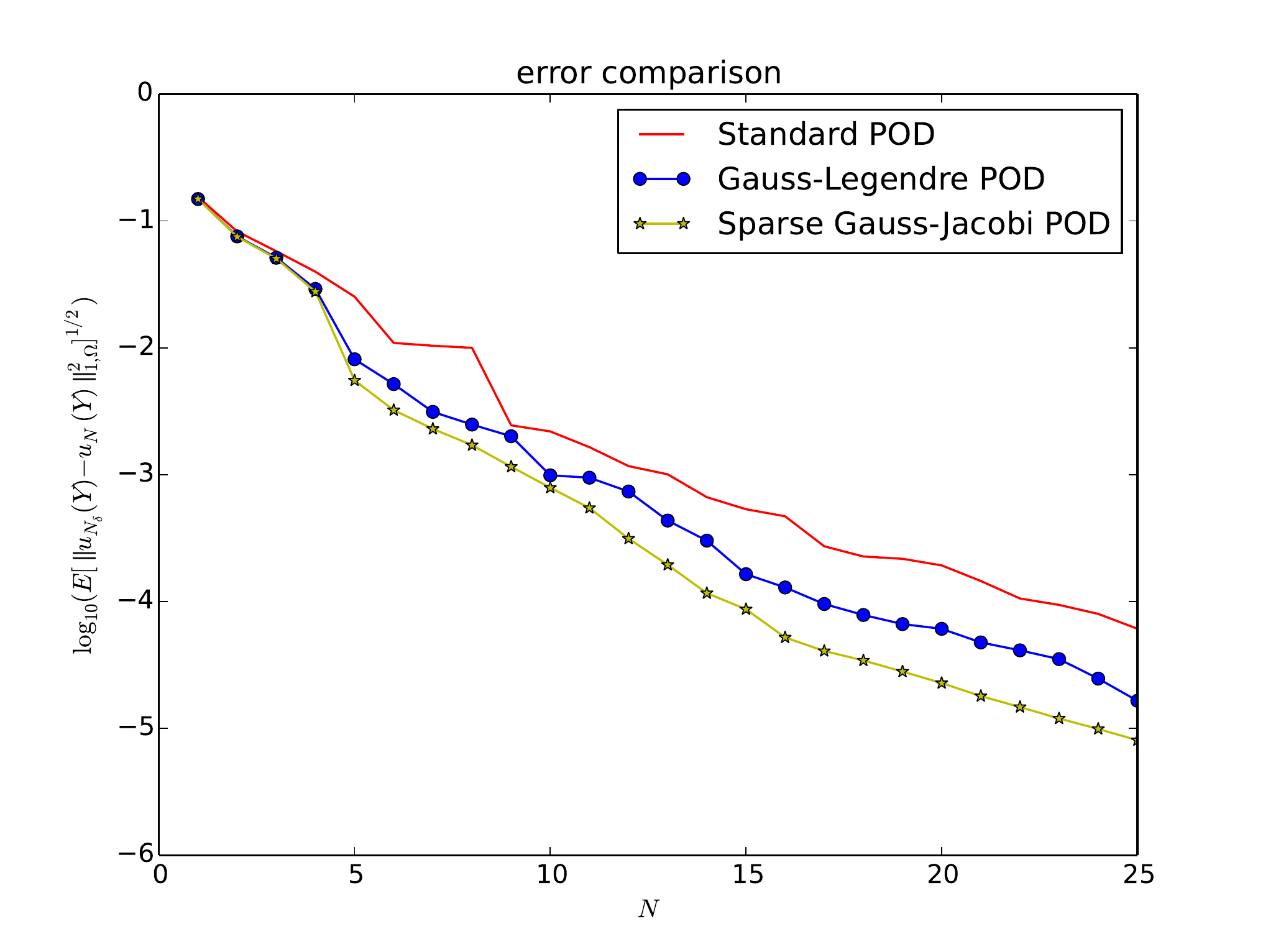}
\end{minipage}
\hspace{4mm}
\begin{minipage}[c]{.47\textwidth}
\includegraphics[width=\textwidth,
keepaspectratio]{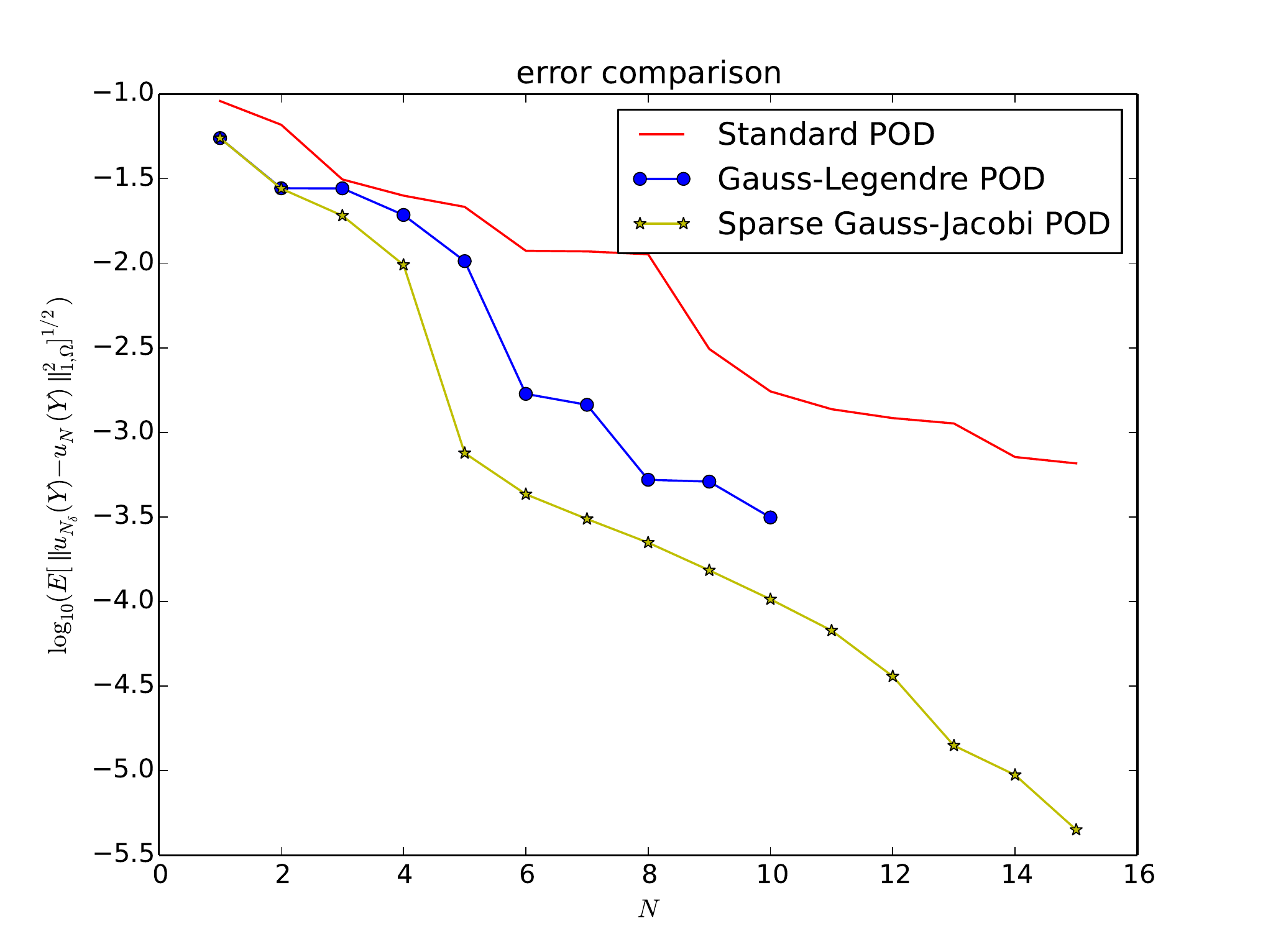}
\end{minipage}
\caption{Error analysis for standard, Gauss-Legendre and sparse Gauss-Jacobi POD algorithms, for the cases $(\alpha, \beta)=(10, 10)$ (left) and $(\alpha, \beta)=(75, 75)$ (right). \label{figure:6E_qPOD}}
\end{figure}

\begin{figure}
\centering
\begin{minipage}[c]{.47\textwidth}
\includegraphics[width=\textwidth,
keepaspectratio]{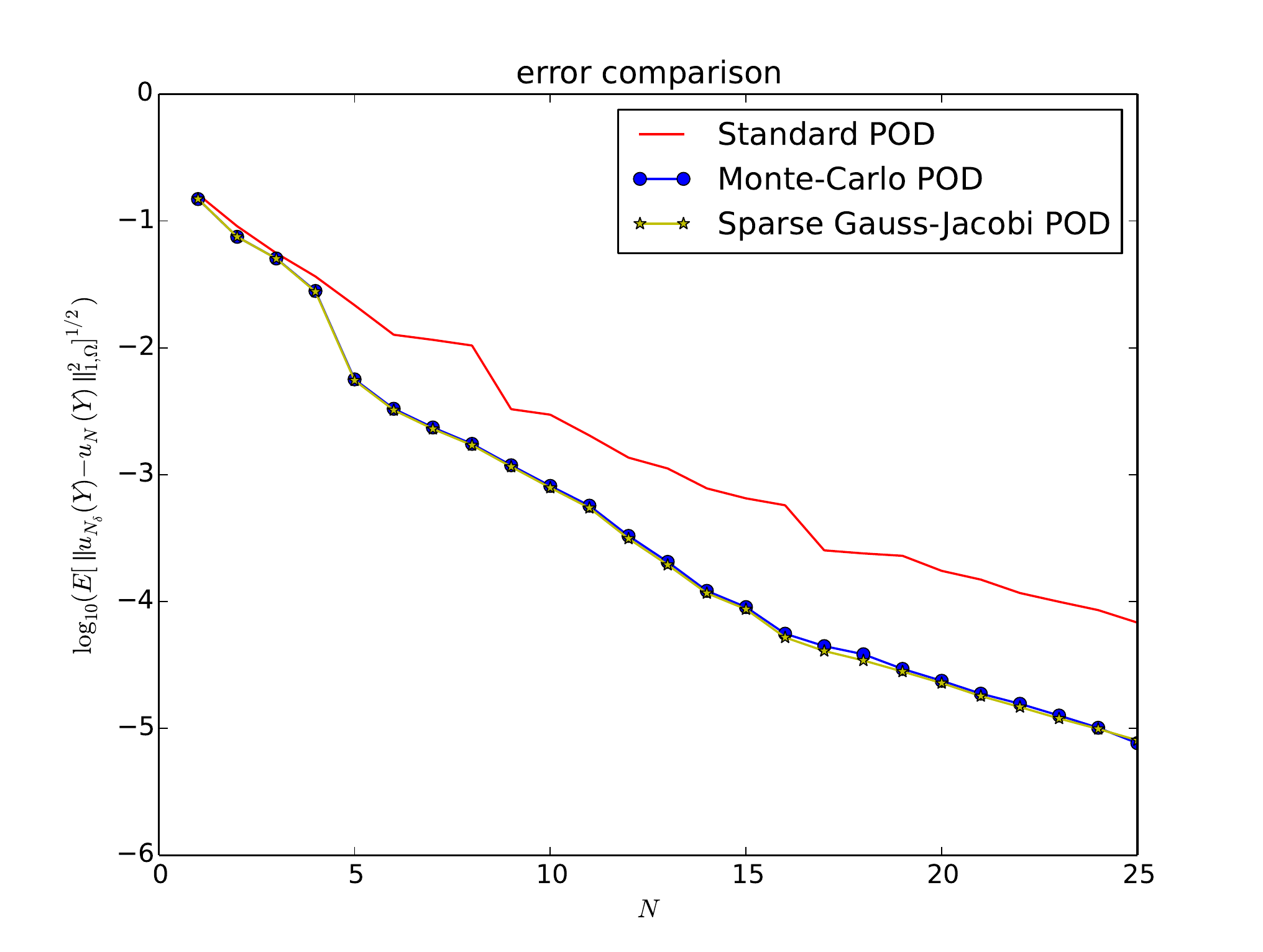}
\end{minipage}
\hspace{4mm}
\begin{minipage}[c]{.47\textwidth}
\includegraphics[width=\textwidth,
keepaspectratio]{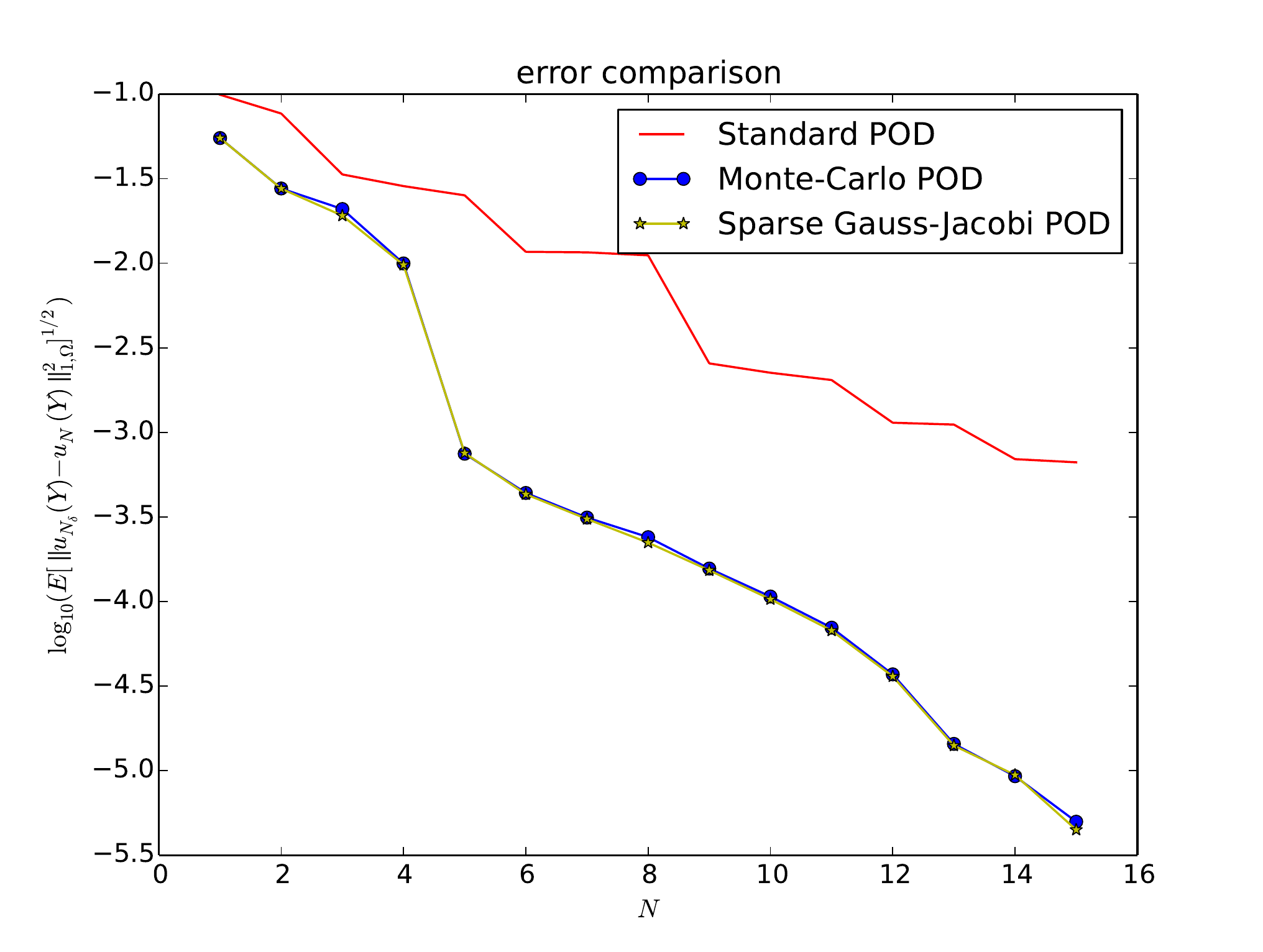}
\end{minipage}
\caption{Error analysis for standard, Monte-Carlo and sparse Gauss-Jacobi POD algorithms,  for the cases $(\alpha, \beta)=(10, 10)$ (left) and $(\alpha, \beta)=(75, 75)$ (right). \label{figure:6E_bPOD}}
\end{figure}

Owing to the interpretation of \eqref{pod_deterministic_l2_error_quad_rule} as a quadrature rule, we provide in \figurename~\ref{figure:6E_qPOD}-\ref{figure:6E_bPOD} a comparison to two further possible POD choices, based on tensor Gauss-Legendre and sparse Gauss-Jacobi quadratures, respectively. The cardinality of the training set in the Gauss-Legendre case is 729, while the sparse Gauss-Jacobi entails 389 training points (resulting from the sparsification of a tensor Gauss-Jacobi rule of 729 nodes), as summarized in \tablename~\ref{table:4EB_POD}. The former value was chosen as the minimum cardinality larger than 500 associated to a tensor Gauss-Legendre quadrature formula, in order to provide a fair comparison to the Monte-Carlo cases discussed previously. \figurename~\ref{figure:6E_qPOD} shows indeed similar improvements as in \figurename~\ref{figure:6E_mcPOD} for what concerns the comparison between the standard POD and the weighted approaches. As in the weighted RB method, an inaccurate exploration of the parameter space caused  the Gauss-Legendre POD to generate singular reduced order systems for $(\alpha, \beta) = (75, 75)$ for $N \geq 10$. Finally, \figurename~\ref{figure:6E_bPOD} depicts a similar comparison between Monte-Carlo and sparse Gauss-Jacobi POD. The latter should be preferred because, even though the two weighted quadrature rules perform similarly in terms of error, the sparse Gauss-Jacobi is obtained through a less dense training set $\trainset$, resulting in a less expensive Offline phase (see lines 2-4 of Algorithm \ref{algo:wpod}).

\section{Conclusion}
\label{sec:conclusion}
In this chapter we have presented two weighted reduced order methods for stochastic partial differential equations. A parametrized formulation has been introduced in order to profit by reduction methods based on the projection onto a low dimensional manifold, which has been identified by means of either a Greedy algorithm or a proper orthogonal decomposition, customarily employed in a deterministic setting. Corresponding weighted approaches have been introduced to properly account for the stochastic nature of the problem. The main rationale for the weighting procedure has been emphasized, especially in relation to the minimization of the expectation of the mean square error between the high fidelity and reduced order approximations. Suitable discretizations of such expectation have been discussed, highlighting in particular the role of (possibly sparse) quadrature rules. A numerical test case for a stochastic linear elasticity problem has been presented. The results show that the weighted methods perform significantly better than the standard (deterministic) ones; bigger improvements have been obtained for more concentrated parameter distributions. Moreover, results demonstrate that the choice of a training set $\Xi_t$ which is representative of the distribution is also essential. Further developments of this work concern the investigation of the proposed weighted reduced order methods in stochastic (nonlinear) fluid dynamics problems.

\begin{acknowledgement}
We acknowledge the support by European Union Funding for Research and Innovation - Horizon 2020 Program - in the framework of European Research Council Executive Agency: H2020 ERC Consolidator Grant 2015 AROMA-CFD project 681447 ``Advanced Reduced Order Methods with Applications in Computational Fluid Dynamics''. We also acknowledge the INDAM-GNCS projects ``Metodi numerici avanzati combinati con tecniche di riduzione computazionale per PDEs parametrizzate e applicazioni'' and ``Numerical methods for model order reduction of PDEs''. The computations in this work have been performed with RBniCS \cite{rbnics} library, developed at SISSA mathLab, which is an implementation in FEniCS \cite{logg2012automated} of several reduced order modelling techniques; we acknowledge developers and contributors to both libraries.
\end{acknowledgement}

{
\small
\bibliography{biblio} 
\bibliographystyle{spmpsci}
}
\end{document}